\documentclass[10pt]{article}
\usepackage[utf8]{inputenc}

\usepackage{amsfonts}
\usepackage{amsmath}
\usepackage{amsthm}
\usepackage{graphicx}
\newcommand{\R}{\mathbb{R}}

\newtheorem{thm}{Theorem}

\newcommand{\tr}{\mbox{tr}}
\newcommand{\disc}{\mbox{disc}}
\title{Global bifurcation map of the homogeneus states in the Gray-Scott model}
\author{Joaqu\'{\i}n Delgado\thanks{Departamento de Matemáticas, Universidad Aut\'onoma Metropolitana--Iztapalapa}
\and
Luc\'{\i}a Ivonne Hern\'andez--Mart\'{\i}nez\thanks{Universidad Aut\'onoma de la Ciudad de M\'exico--San Lorenzo Tezonco.},
\and
 Javier P\'erez--L\'opez\thanks{Departamento de Matemáticas, Universidad Aut\'onoma Metropolitana--Iztapalapa}}

\begin{document}
\maketitle

\begin{abstract}
We study the spatially homogeneous time--dependent solutions and their bifurcations  of the Gray-Scott model. We find the global map of bifurcations by a combination of rigorous verification of the existence of  Takens-Bogdanov and a Bautin bifurcations, in the space of two parameters $k$--$F$. With the aid of numerical continuation of local bifurcation curves we give a global description of all the possible bifurcations.
\end{abstract}

\section{Introduction}
The  Gray-Scott model is a system of two reaction--difussion partial differential equations depending on two parameters that describes an auto--catalytic chemical reaction similar to  the famous Belouzov--Zhabotinzky. Its numerical study, initiated by Pearson \cite{Pearson}  has been  continued extensively by Munafo \cite{Munafo} among other authors. Some investigations include extensions of the model such as  as a stochastic component  \cite{Wang}, or pose it as   a control problem \cite{kyrychko}. In \cite{Mazin} the homogeneous states are described succinctly and some other patterns bifurcating from them are shown. In most of these papers the homogeneous state defined by the trivial critical point $u=1$, $v=0$ is perturbed by changing the values of $u$ and $v$ in a small central region and then perturbing randomly in order to break any a priori symmetry. The numerical experiments by these an other authors show  that a rich variety of patterns develop when   a particular state is perturbed. A zoo of patterns such as spots, stripes, rings etc. are common.  In \cite{Xmorphia} a nice map of  numerical explorations in the parameter space $k$--$F$ is presented on-line. 

Instabilities from the homogeneous states arise by different recognized mechanisms such as Turing and Hopf. The common floor in these approaches relies in a full understanding of the spatially homogeneous states which are taken as initial solution to be perturbed. These are solutions of (\ref{Gray-ScottEDP}) that do not depend on $(x,y)$ but solely on time $t$, and therefore are solutions of the ODE obtained by discarding the diffusion terms. Several authors mention  the importance of  the locus in parameter spaces of Hopf and saddle--node bifurcations and particularly of a Takens--Bogdanov bifurcation at the particular value $(k,F)=(1/16,1/16)$ and the change in sign of the first Liapunov coefficient, leading to unstable or stable limit cycles por an approximate value  of $k=0.035$.

Nevertheless a review of the available literature is lacking of formal proofs on the existence of a Takens--Bogdanov (BT)  bifurcation. Also the particular value were the first Liapunov coefficient vanishes suggests the possibility of a Generalized Hopf bifurcation. 
This is the main contribution of the present paper which is two--fold. We provide rigorous proofs of the existence of a BT bifurcation and of a Bautin bifurcation; in the second we provide the exact value $(k,F)= \left(9/256,3/256\right)$ for which this type of bifurcation occurs also. Secondly, knowing the known local bifurcation diagrams  of the respective normal form, and numerical continuation, we provide a global map of bifurcation. In particular we prove the existence of two limit cycles. The bifurcation curves then divide the parameter space into regions, were qualitatively the phase portrait is described.

The rest of the paper is organized as follows: In Section 2 we describe the system to be studied in the paper and the critical points. In Section 2 we give a geometric description of the bifurcation sets, particular of the saddle--node bifurcations. In Section 3 we provide analytical formulas for the  critical points. In section 4 we analyze the linear stability. In section 5 and 6  we state Theorems 1 and 2, respectively  on the BT and Bautin bifurcations. There we give the main ideas of the proofs and postponed the technical details to the Appendix.  In Section 7 we proceed to describe the global map of bifurcation which can be safely summarized schematically in Figure~11. Then we present the phase portraits for each of the regions determined by the bifurcation curves.

\section{The spatially homogeneous states}
The Gray-Scott model is a system of reaction-diffusion PDEs representing the inflow-outflow, reaction and diffusion of two chemical components with concentrations $u$, $v$, given by
\begin{eqnarray*}
u_t &=& D_u \nabla^2 u- uv^2 + F(1-u),\nonumber\\
v_t &=& D_v \nabla^2 v + uv^2 - (F+k) v.\label{Gray-ScottEDP} 
\end{eqnarray*}
Here $u$, $v$ are function depending on the time and space $(t,x,y)$, $t\geq 0$, and   $(x,y)$ varying in the rectangular domain
$\Omega=[0,1]^2.$ The diffusion coefficients $D_u$, $D_v$ are constants, and the  positives parameters $F$ and $k$, represent the inflow of substance $u$ and the outflow of substance $v$. Usually Neumann or periodic  boundary conditions are used but these will fixed later.

On a first stage we study the spatially homogeneous states that depend only on time. This yields the system of ODE
\begin{eqnarray}\label{system}
u' &=& - uv^2 + F(1-u),\\
v' &=&  uv^2 - (F+k) v.\nonumber
\end{eqnarray}

We can see that the first quadrant $u>0$, $v>0$ is invariant, and therefore we will consider the system defined for $u\geq 0$, $v\geq0$.
In what follows we write this system as
$$
\dot{p}=f(p,\alpha),
$$
where $p= (u,v)^T$, $\alpha=(k,F)^T$ is the vector of parameters. We will omit the reference to the vector of parameters when its mention is not specially important and use the shorthand $\dot{p}=f(p)$.

\section{The surface of critical points, singular and  bifurcation sets}
\begin{figure}
\centering
\includegraphics[scale=.3]{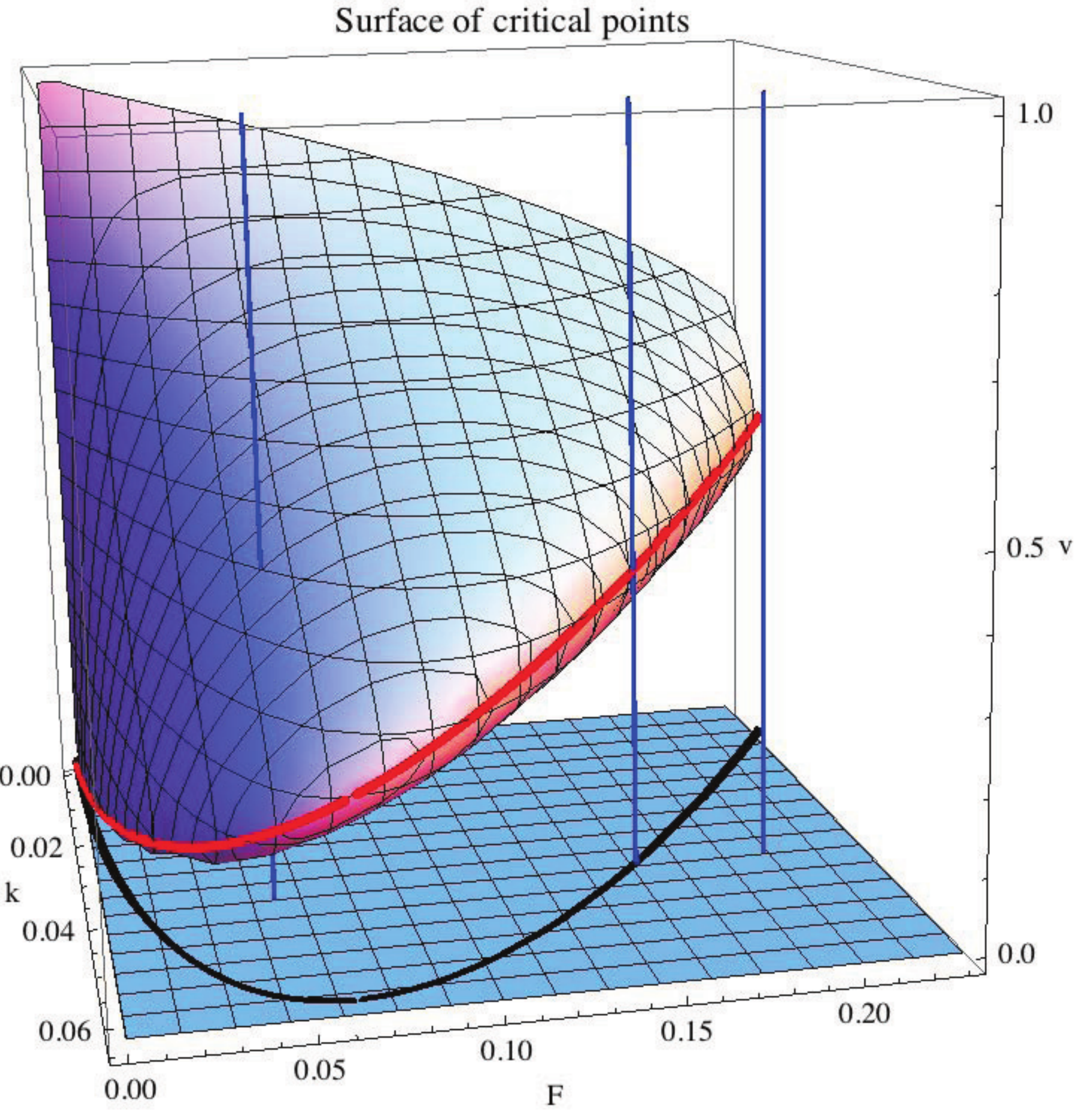}
\caption{The set of critical points. The trivial component $v=0$ yields the trivial point $(u,v)=(1,0)$. The singular set is the red curve on the non-trivial component $\Sigma$ and its projection, the bifurcation set is the black curve in parameter space $(k,F)$. For each value of $(k,F)$ there are 1, 2 or 3 critical points, depicted in the figure by vertical blue lines.\label{surface}}
\end{figure}

Geometrically the number of critical points can be described as follows: Let $f=(f_1,f_2)$ denote the vector field defining the system  (\ref{system}), then as $f_i$ are functions of  $(k,F,u,v)$ and polynomial in $(u,v)$ we can take the resultant with respect to the variable $v$,
$$
R= \mbox{Res}[f_1,f_2,u]=v \left(F (F + k) - F v + (F + k) v^2\right)
$$
which is a polynomial in $v$ with coefficients depending on $(k,F)$.  Therefore, the surface $R(k,F,v)=0$ in the space $k$--$F$--$v$ describes the number of critical points as we vary the parameters $(k,F)$. The component $v=0$  yields the trivial critical point $(u,v)= (1,0)$. Let  the non trivial component be defined by
$$
\Sigma =\{(k,F,v)\mid G(k,F,v)=F (F + k) - F v + (F + k) v^2=0\}.
$$
and $\pi\colon\R^3\to\R^2$, $(k,F,v)\mapsto(k,F)$ the projection to  parameter space. The singular set is the curve on $\Sigma$ where the tangent plane to $\Sigma$ is parallel to the $v$-axis or
$$
0= \nabla G(k,F,v)\cdot (0,0,1)=-F + 2 (F + k) v
$$
Thus the singular set is the intersection of the surfaces $\Sigma$ and  $-F + 2 (F + k) v=0$, and its projection into the parameter space $(k,F)$ is the bifurcation set. By eliminating $v$ using the resultant with respect to $v$ set we obtain  $4(F+k)^2-F=0$ or equivalently $\Delta=0$, where $\Delta$ is the discriminant defined below (\ref{discriminant}). Figure~\ref{surface}  describes the set of critical points and the singular  and bifurcation sets in $\Sigma$. The red curve is the singular set and divides $\Sigma$ in two components.  The critical point on the upper component will be denoted by $p_{\pm}$ and the critical point on the lower component will be denoted by $p_{\mp}$. The trivial critical point will be denoted by $p_{0}$. The bifurcation set is the black curve in the plane of parameters $k$--$F$.

\section{Equilibrium points and their stability}
As it was described  in the previous section the trivial critical point $(u,v)=(1,0)$ exist for all values of the parameters $(k,F)$.

Supposing now $v\neq0$, we equate both equations in  (\ref{system}) to zero, and cancel out the factor $v$ in the second, then sum  the right hand sides to get the equivalent system for the critical points:
\begin{eqnarray}
u+\gamma v&=&1,\label{line}\\
uv&=&F+k. \label{hyper}
\end{eqnarray}
where $\gamma=(F+k)/F$. The first equation
says geometrically  that the critical points $p_{\mp}$ and $p_{\pm}$ lay along a straight line
and at the intersection of the hyperbola (\ref{hyper}) with this line. The first condition in turn implies, since $u,v$ are non--negative, that while $0<u<1$ then $0< v < \gamma^{-1}$. Substituting $u$ from (\ref{line}) into (\ref{hyper}) yields the quadratic equation
\begin{equation}\label{quadratic}
v^2-\frac{1}{\gamma}v+F=0
\end{equation}
from which we get
$$
v_{\pm}=\frac{1}{2\gamma}(1\pm\sqrt{\Delta})
$$
where
\begin{equation}\label{discriminant}
 \Delta=1-4F\gamma^2.
\end{equation}
Substituting these values of $v_{\pm}$  into (\ref{line}) we get the corresponding values for $u$.

In summary, we have
if $\Delta>0$ (it is necessarily  that $0<F<1/4$), there are two critical points given by
\begin{eqnarray}
p_{\mp}=(u_{-},v_{+})&=&\left(\dfrac{1}{2}\left(1-\sqrt{\Delta}\right),\dfrac{1}{2\gamma}\left(1+\sqrt{\Delta}\right)\right),\label{cp1}\\
p_{\pm}=(u_{+},v_{-})&=&\left(\dfrac{1}{2}\left(1+\sqrt{\Delta}\right),\dfrac{1}{2\gamma}\left(1-\sqrt{\Delta}\right)\right).\label{cp2}
\end{eqnarray}
If $\Delta=0$ the two critical points coincide,
$$
p_{\pm}=p_{\mp}= \left(\frac{1}{2},\frac{1}{2\gamma}\right)
$$
and if $\Delta<0$ there exists no additional critical points.

Figure~\ref{fig1} shows the curve $\Delta=0$ (green curve),  dividing the plane positive quadrant $k$--$F$ in two components. The bounded component corresponds to $\Delta>0$, where there exists three critical points, on the unbounded component $\Delta<0$ and only the trivial critical point exist. On the bifurcation set $\Delta=0$, two critical points coincide and we will show that they are saddle--nodes.

\section{Linear stability}
The Jacobian matrix of (\ref{system}) at any critical point $p=(u,v)$ is
\begin{equation}\label{linearization}
Df(p)= \left(
\begin{array}{cc}
-(F+v^2) & -2uv\\
v^2 & -(F+k)+2uv
\end{array}\right)
\end{equation}

The Jacobian matrix (\ref{linearization}) at the critical point $p_0=(1,0)$ yields
$$
Df(p_{0})= \left(
\begin{array}{cc}
-F & 0\\
0 & -(F+k)
\end{array}\right)
$$
therefore the eigenvalues at  are $-F$ and $-(F+k)$ and  trivial critical point is always asymptotically stable.

From (\ref{cp1}), (\ref{cp2}) it follows that for the critical points $p_{\mp}$ and $p_{\pm}$
$$
u_{-}v_{+}=u_{+}v_{-} = \frac{1-\Delta}{\gamma}=F+k,
$$
therefore the linearization at any of these point is
$$
Df(p)= \left(
\begin{array}{cc}
-(F+v^2) & -2(F+k)\\
v^2 & (F+k)
\end{array}\right).
$$
In particular the trace and determinant are
\begin{eqnarray}
\det(Df(p)) &=& (v^2-F)(F+k),\\
\mbox{tr}(Df(p)) &=& k-v^2\label{traza},
\end{eqnarray}
where $v$ is evaluated at  the corresponding  the critical point.

\subsection{The saddle--node curve}
Substituting $v$ for each of the critical point (\ref{cp1}) and (\ref{cp2}) we get for the determinant
\begin{eqnarray}
\det(Df(p_{\pm}))&=&\dfrac{F \left(\Delta -\sqrt{\Delta}\right)}{2\gamma},\label{ec:detp1}\\
\det(Df(p_{\mp}))&=&\dfrac{F \left(\Delta +\sqrt{\Delta}\right)}{2\gamma},\label{ec:detp2}
\end{eqnarray}

Since $k,F$ are positive it follows that $\Delta<1$, therefore we see that $\Delta=0$ if and only if $p_{\pm}=p_{\mp}$ and $\det(Df(p_{\pm}))=\det(Df(p_{\mp}))=0$, thus the bifurcation curve $\Delta=0$ is a curve of saddle--nodes where at least one eigenvalues vanishes. Solving the equation $\Delta=0$ in terms of $F$ yields the two branches
\begin{equation}
F= \frac{1}{8} (1 - 8 k) \pm \sqrt{1 - 16 k},\label{branches}\\
\end{equation}
which are shown in Figure~\ref{fig1}.

\subsection{Stability of $p_{\pm}$ and $p_{\mp}$}
Computing the trace according to (\ref{traza}) yields
\begin{equation}
\text{tr}(Df(p))=\frac{2 \gamma ^3 F-1\pm\sqrt{\Delta}}{2 \gamma ^2}.\label{ec:trp1}\\
\end{equation}
where the upper signs refers to the critical point $p=p_{\pm}$ and the lower sign to $p=p_{\mp}$.

\begin{figure}
\centering
\includegraphics[width=0.7\textwidth]{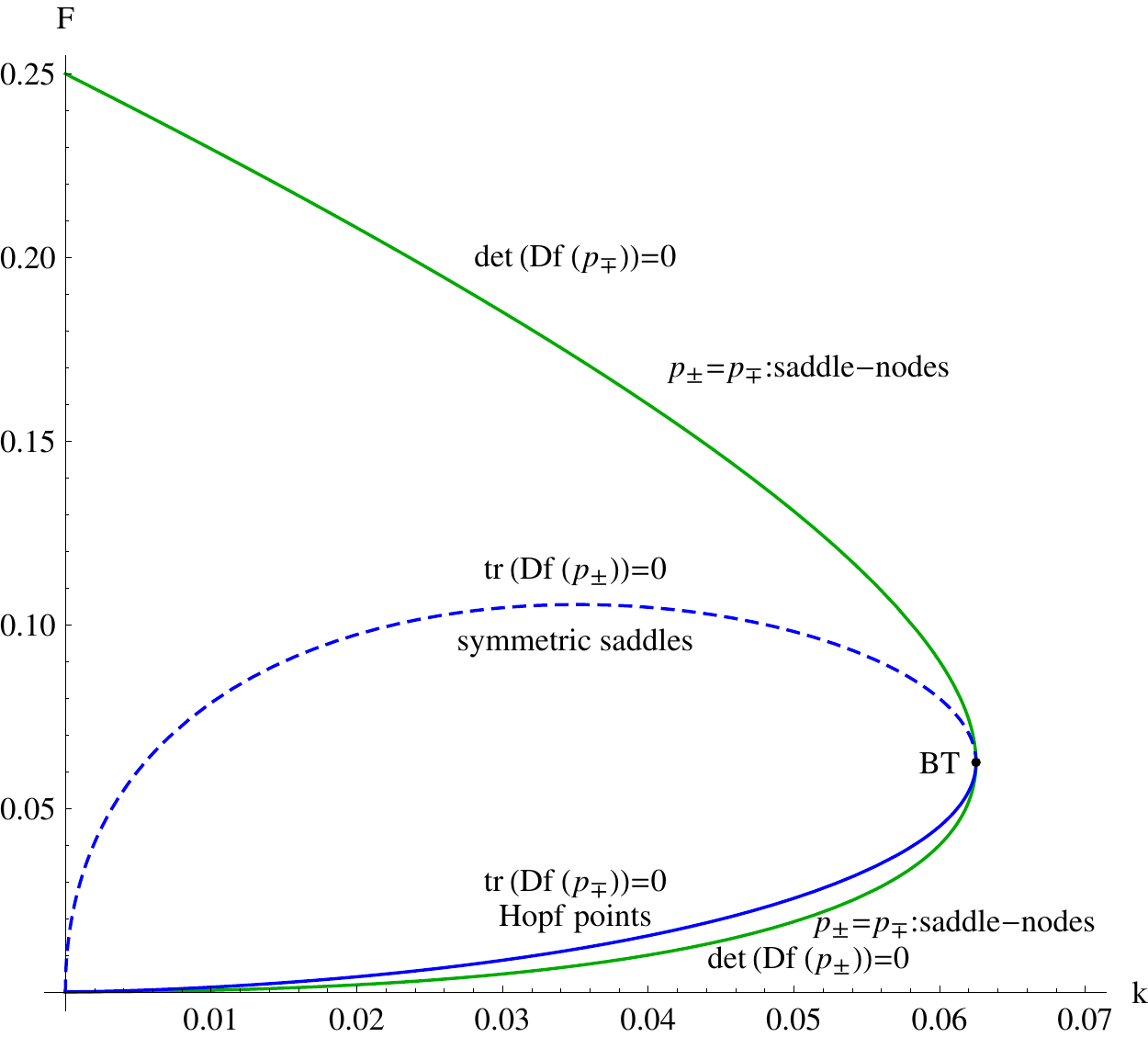}
\caption{The curves of determinant and trace equal to zero. \label{fig1}}
\end{figure}

Solving for $F$ we see that the curves  $\tr(Df(p))=0$ can be parameterized as
\begin{equation}
F = \frac{1}{2}\left(\sqrt{k}-2k \pm \sqrt{k-4k\sqrt{k}}\right),\label{curve:symmetric_saddle}
\end{equation}
with the same convention of signs as before.
These curves are shown as green curves  in Figure~\ref{fig1}. There is a particular critical point point where both $\tr(Df(p))$ and $\det(Df(p))$ vanish. We will verify (see Section~\ref{BT-bifurcation}) that it is in fact a Bogdanov-Takens point denoted by BT in Figure~\ref{fig1}. Solving for this special point we get the values of the parameters and coordinates of the critical point as
\begin{equation}\label{BT-point}
BT=(k,F)=\left(\frac{1}{16},\frac{1}{16}\right),\qquad
p_{BT}= (u,v)= \left(\frac{1}{2},\frac{1}{4}\right).
\end{equation}

\begin{figure}
\centering
\includegraphics[width=0.7\textwidth]{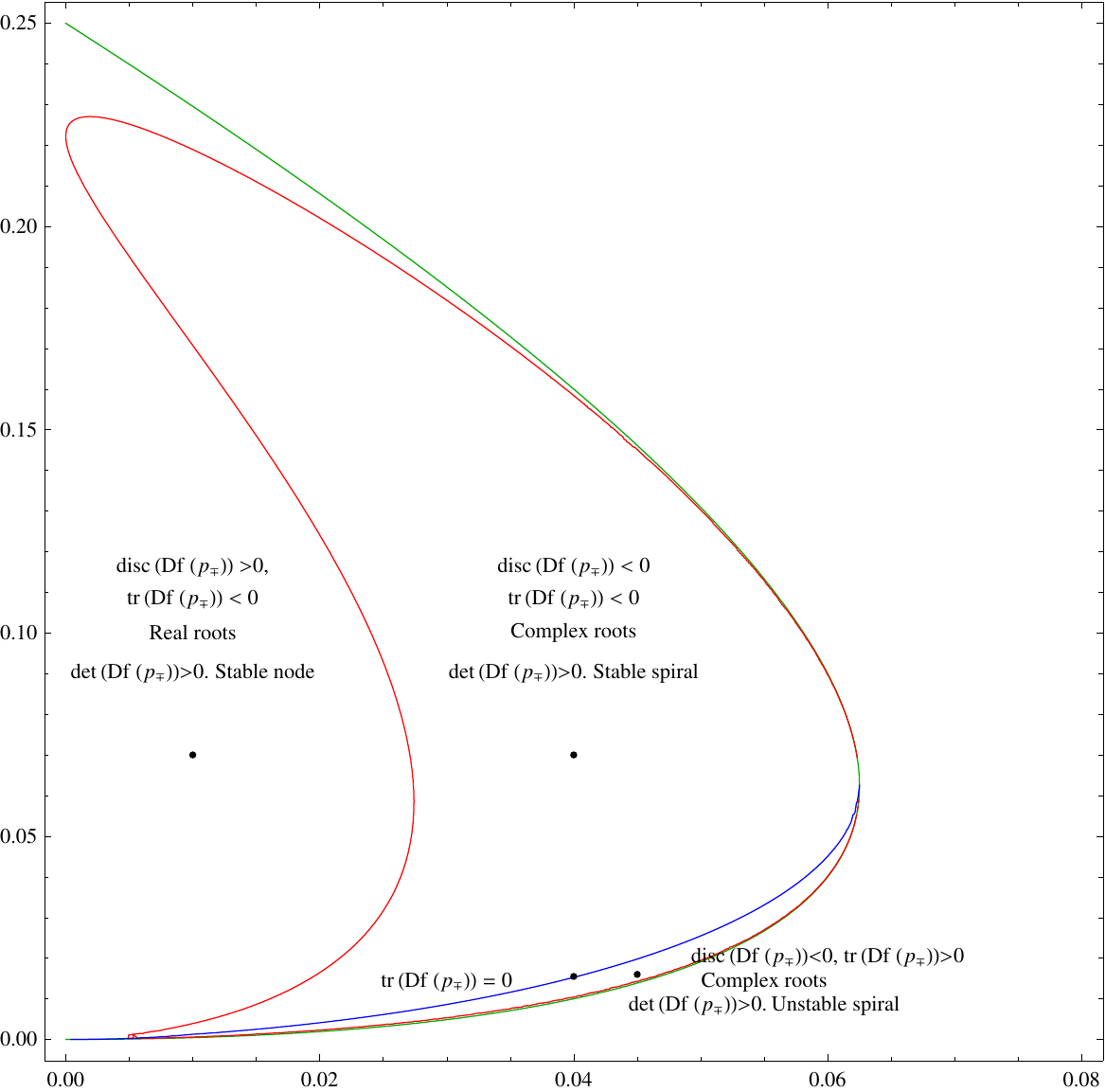}
\caption{
The curves $\text{tr}(Df(p_{\mp}))=0$ (blue),  $\disc(Df(p_{\mp}))=0$ (red)  for the critical point $p_{\mp}$.
\label{Delta2}}
\includegraphics[width=0.7\textwidth]{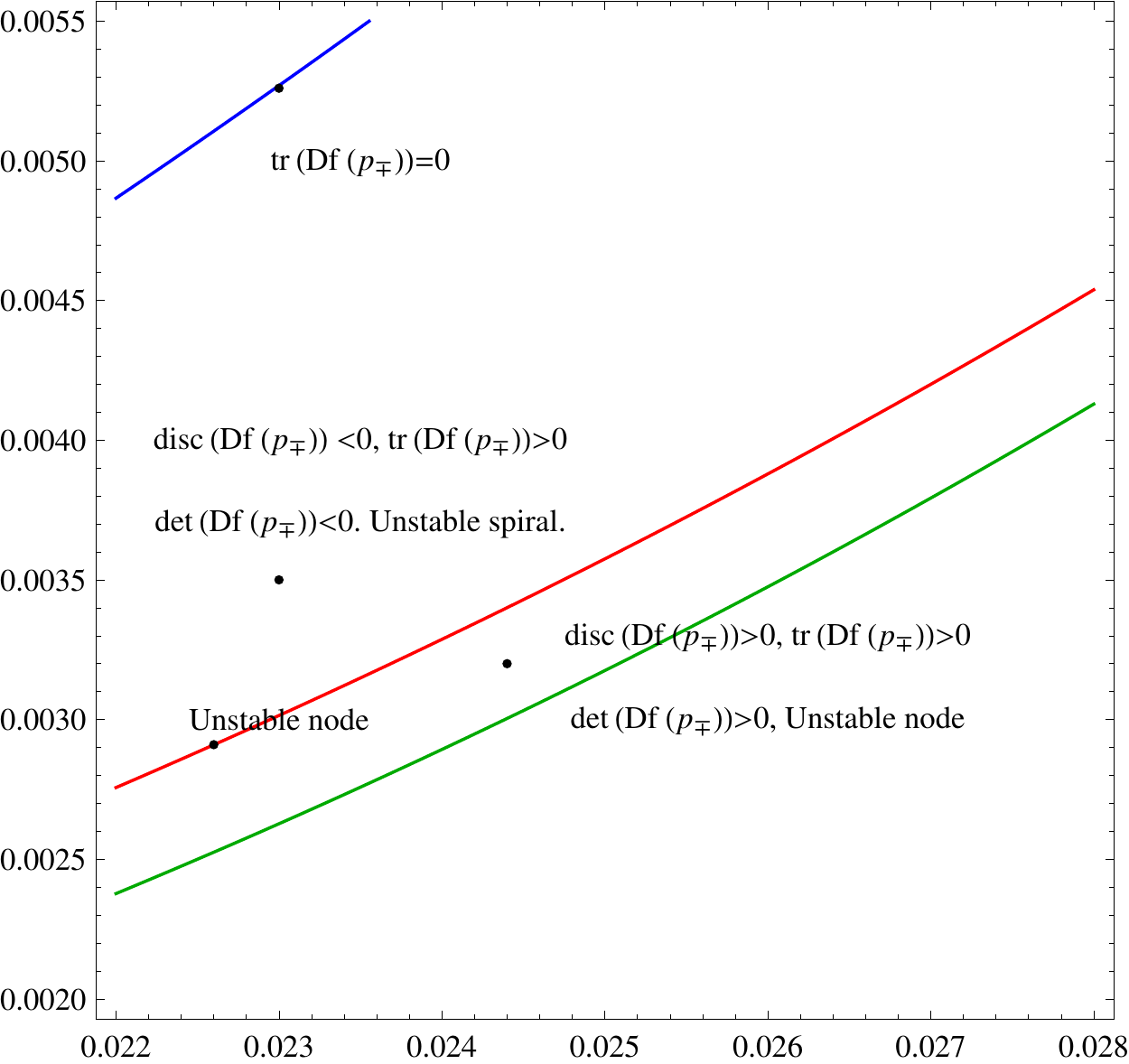}
\caption{
Zoom of figure \ref{Delta2} between the Hopf curve and the discriminant curve.
\label{Delta2_zoom}}
\end{figure}

 For $p_{\mp}$, it follows from (\ref{ec:detp2}) that $\det( Df(p_{\mp}))>0$. Here we have the posibility of complex eigenvalues, in particular if $\tr(Df(p_{\mp}))=0$, the real parts are zero. Solving (\ref{ec:trp2}) for $F$ gives the curve (\ref{curve:symmetric_saddle}) with the lower sign (corresponding to this case $p_{\mp}$). Thus
 \begin{equation}
 F = \frac{1}{2}\left(\sqrt{k}-2k - \sqrt{k-4k\sqrt{k}}\right)
 \end{equation}
parameterizes the curve of Hopf points and correspondingly,
 \begin{equation}
 F = \frac{1}{2}\left(\sqrt{k}-2k + \sqrt{k-4k\sqrt{k}}\right)
 \end{equation}
parameterizes the curve of symmetric saddles. These curves are shown in Figure~\ref{Delta2} as blue continuous and dashed curve respectively.

\medskip
The discriminant $\disc(Df(p_{\mp}))  = \tr(Df(p_{\mp}))^2-4\det(Df(p_{\mp})),$
changes sign along the curve $\disc(Df(p_{\mp}))=0$ shown in red in Figure~\ref{Delta2}.
Inside the bounded region $\disc(Df(p_{\mp}))<0$ and we have complex roots. This curve is tangent to the saddle--node curve at  the Takens--Bogdanov point   and lies between the Hopf and the saddle--node curves in its lower part as seen in more detail in  Figure~\ref{Delta2_zoom}.

Let us now describe how the stability of the critical point $p_{\mp}$  changes as we follow an arbitrary line $k=const$ decreasing the value of $F$. Above the Hopf curve $\tr(Df(p_{\mp}))=0$ (shown in blue in Figure~\ref{Delta2}), the critical point $p_{\mp}$ is stable. Inside the region $\disc(Df(p_{\mp}))<0$ it is stable spiral; outside the region $\disc(Df(p_{\mp}))>0$ its a stable node.  When the point approaches  the Hopf curve, the eigenvalues are complex and the critical point $p_{\mp}$ becomes a stable spiral. When the point crosses the Hopf curve $\tr(Df(p_{\mp}))=0$  a limit cycle is formed through a Hopf bifurcation (we will determine  later its stability through the computation of the first Liapunov coefficient). Below the Hopf curve and above in the curve $\disc(Df(p_{\mp}))=0$ the point $p_{\mp}$ becomes an unstable spiral. Decreasing the value of $F$, the critical point encounters the discriminant curve  $\disc(Df(p_{\mp}))=0 $ and then the critical point becomes an unstable node.  Finally, the point $p_{\mp}$ collides with the point $p_{\pm}$ along the saddle--node curve $\det(Df(p_{\mp}))=0 $ in its lower component.

\section{Takens-Bogdanov bifurcation\label{BT-bifurcation}}
Several authors \cite{kyrychko}, \cite{Pearson}, \cite{Wang}, among others report a  Takens-Bogdanov bifurcation point  $BT= (k,F)=\left(1/16,1/16\right)$.
Here we give a rigorous proof of the existence of a Takens--Bogdanov bifurcation by an explicit verification of the non--degeneracy conditions as in \cite{Kuznetsov}. This has the implication that locally the system is topologically equivalent to one of its known normal forms.

Recall that the curve of saddle nodes is defined by the equation $\Delta(k,F)=0$.

\begin{thm}\label{tTB}
System (\ref{system}) undergoes a Takens--Bogdanov bifurcation at $BT\equiv(k,F)=\left(1/16,1/16\right)$. As a consequence, in the parameter space $k$--$F$, there exists local  curve of Hopf bifurcations  and a local  curve of  homoclinic bifurcations. These curves are tangent to the saddle node curve defined by the equation $\Delta(k,F)=0$ at $BT$.
\end{thm}

For $(k,F)=\left(1/16,1/16\right)$ the critical point is  $p_{\pm}=p_{\mp}=\left(1/2,1/4\right)$. Let us perform a shift of the origin of coordinates and parameters,
$$\begin{array}{rclcrcl}
\alpha_1&=&F-1/16,&\qquad x_1&=&u-1/2,\\
\alpha_2&=&k-1/16,&\qquad x_2&=&v-1/4,
\end{array}$$
then the system (\ref{system}) becomes
\begin{eqnarray}
\dot x_1&=&-\left(\dfrac{(\alpha_2+\frac{1}{16})^2}{4(\alpha_1+\alpha_2+\frac{1}{8})^2}+\alpha_2+\frac{1}{16}\right)x_1-\dfrac{\alpha_2+\frac{1}{16}}{2(\alpha_2+\alpha_1+1/8)}x_2\nonumber\\
&-&\dfrac{\alpha_2+\frac{1}{16}}{\alpha_2+\alpha_1+\frac{1}{8}}x_1x_2-\dfrac{1}{2}{x_2}^2-x_1{x_2}^2,\nonumber \\
\dot x_2&=&\dfrac{(\alpha_2+\frac{1}{16})^2}{4(\alpha_2+\alpha_1+\frac{1}{8})^2}x_1+\left(\dfrac{\alpha_2+\frac{1}{16}-2(\alpha_2+\alpha_1+\frac{1}{8})^2}{2(\alpha_2+\alpha_1+\frac{1}{8})}\right)x_2\label{df}\\
&+&\dfrac{1}{2}x_2^2+\dfrac{\alpha_2+\frac{1}{16}}{\alpha_2+\alpha_1+\frac{1}{8}}x_1x_2+x_1x_2^2,\nonumber
\end{eqnarray}
Let us denote this system as $\dot x=f(x,\alpha)$. The idea of the proof is to compute the right and left eigenvectors and generalized eigenvectors of the linear part $Df(0,0)$  and expand up to third order in these coordinates,
then we verify the non degeneracy conditions \cite[p. 321]{Kuznetsov}. Details of the proof are given the Appendix \ref{TB}. In particular the local bifurcation diagram is topologically equivalent to the normal form
\begin{eqnarray*}
\dot{\eta}_1 &=&  \eta_2,\\
\dot{\eta}_2 &=&  \beta_1 + \beta_2 \eta_1 + \eta_2^2 + s \eta_1 \eta_2,
\end{eqnarray*}
where $s=\mbox{sign}(b_{20}(0)a_{20}(0)+b_{11}(0))$,  and the coefficients are given Appendix (\ref{TB}). In our case $s=-1$, so the qualitative diagram is the same as \cite[p. 322]{Kuznetsov}.

\section{Bautin bifurcation}
Several authors \cite{Mazin}, \cite{Pearson}, \cite{Wang},  have pointed out
the stability of limit cycles bifurcating from the Hopf curve, stable for $k<0.035$ and unstable for $k>0.035$, this value has been reported numerically. We give a rigorous proof that for
$(k,F)= \left(9/256,3/256\right)$ a Bautin bifurcation takes place. In particular
this proves the stability of limit cycles as described in the aforementioned references. Moreover, the existence of this bifurcation guarantees the existence
of a curve of limit points of cycles, and the coexistence of stable and unstable limit cycles as described in the local bifurcation diagram of the Bautin bifurcation \cite[p. 313]{Kuznetsov}

\begin{thm}\label{Bautin}
System (\ref{system}) undergoes a Bautin bifurcation at $GH\equiv(k,F)=\left(9/256,3/256\right)$. As a consequence, in the parameter space $k$--$F$, there exists a local  curve of Hopf bifurcations  and a curve of limit point of cycles and an open region where two limit cycles with different stability coexist.
\end{thm}

We give the main idea of the proof. Detailed computations are given in the
Appendix \ref{B}.

The main difficult is to find a manageable expression for the first Liapunov coefficient $\ell_1$ as function of the parameters and then look for its zeros. There are several formulas for $\ell_1$ but most of them requires to put in normal form the linear part, which is cumbersome in our case since the system depends on parameters (see the discussion in \cite{Kuehn}, but notice a couple of typos in formula (12)). The most convenient formula for us is the one in \cite[p. 211]{Chow-Li-Wang}, that we denote by $\ell_1^{CLW}$, which is fully developed in \cite{Marsden}, that does not require the linear normal form. Then we verify the non--degeneracy hypotheses of the Bautin bifurcation as in \cite[p.311]{Kuznetsov} using Kuznetsov formula $\ell_1^{Kuz}$.

For the first part,  we perform a change parameters $(k,F)\mapsto (k,\nu)$ such that $\nu=0$ is the Hopf curve, then we  use expression for the $\ell_1^{CLW}(k,\nu)$, set $\nu=0$ and solve the equation  $\ell^{CLW}(k,0)=0$, we  get
$$
GH=\left(\frac{9}{256},\frac{3}{256}\right), \qquad p_{\mp}=\left(\frac{1}{4},\frac{3}{16}\right).
$$

See details in appendix \ref{B}.

\subsection{Dynamics in a neighborhood of the Bautin point}
\begin{figure}\centering
\includegraphics[width=0.7\textwidth]{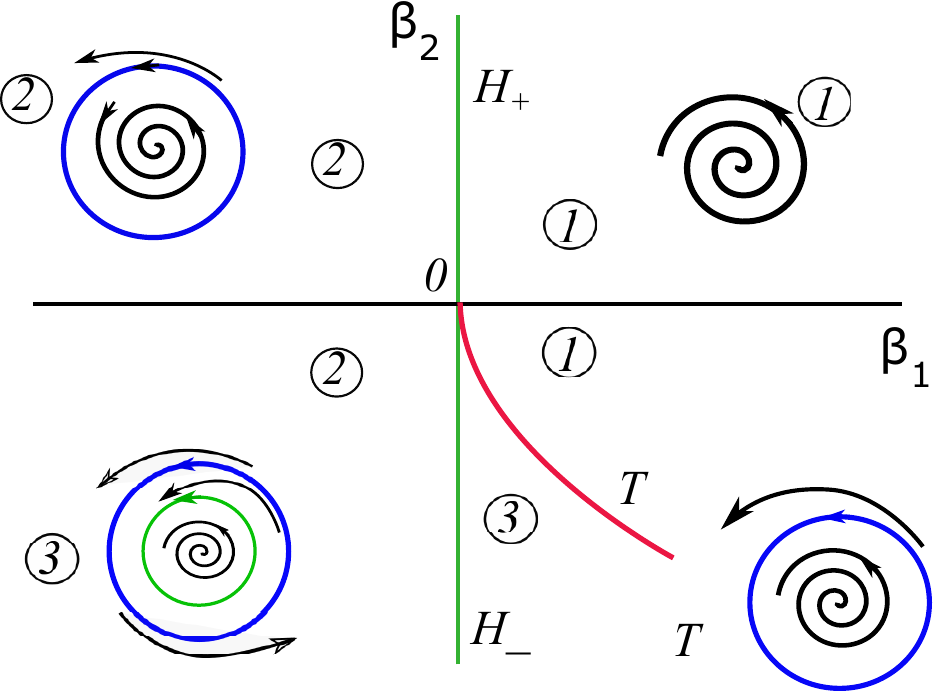}
\caption{Bautin bifurcation for $\ell_{2}$ positive.\label{Bautin-bif}}
\end{figure}

Let us briefly describe the dynamics of the homogenous states in a neighborhood of the Bautin point GH. In Figure \ref{Bautin-bif} there is shown
qualitatively a small neighborhood of  $GH$ in
the space of parameters $k$--$F$.
The local normal form of the system in polar coordinates has the expression
\begin{eqnarray*}
\dot \rho &=& \rho(\beta_1 + \beta_2\rho^2 + \rho^4),\\
\dot \phi &=& 1.
\end{eqnarray*}

The  $\beta_2$ axis is a parametrization of the Hopf curve with the positive
$\beta_2$--axis corresponding to the right of $GH$, the negative $\beta_2$--axis to the left.
In this coordinates, $GH=(0,0)$ and the first Liapunov coefficient is precisely $\beta_2$.
The  curve $T = \{(\beta_1,\beta_2): \beta_2^{2}- 4\beta_1 = 0,\quad \beta_2<0\}$ of limit points of cycles,
according to Theorem~\ref{Bautin}, together with the negative $\beta_2$--axis define a wedge region 3.
Region 2 is defined as the semiplane $\beta_1<0$, and Region 1 is the complement.
Taking a closed path starting in Region 1,
where the system has a single unstable equilibrium and no cycles. Crossing the Hopf bifurcation
boundary $H_{+}$ from region 1 to region 2  the equilibrium point becomes stable, and  a unique
an unstable limit cycle appears, witch survives when we enter region 3.
Crossing the Hopf boundary $H_{-}$ creates an extra stable cycle inside the first one, while the equilibrium point is unstable. Two cycles
of opposite stability exist inside region 3 and disappear at the curve $T$,
in which the point of equilibrium and the cycle are unstable, thus completing the circle.

\begin{figure}
\centering
\includegraphics[scale=.5]{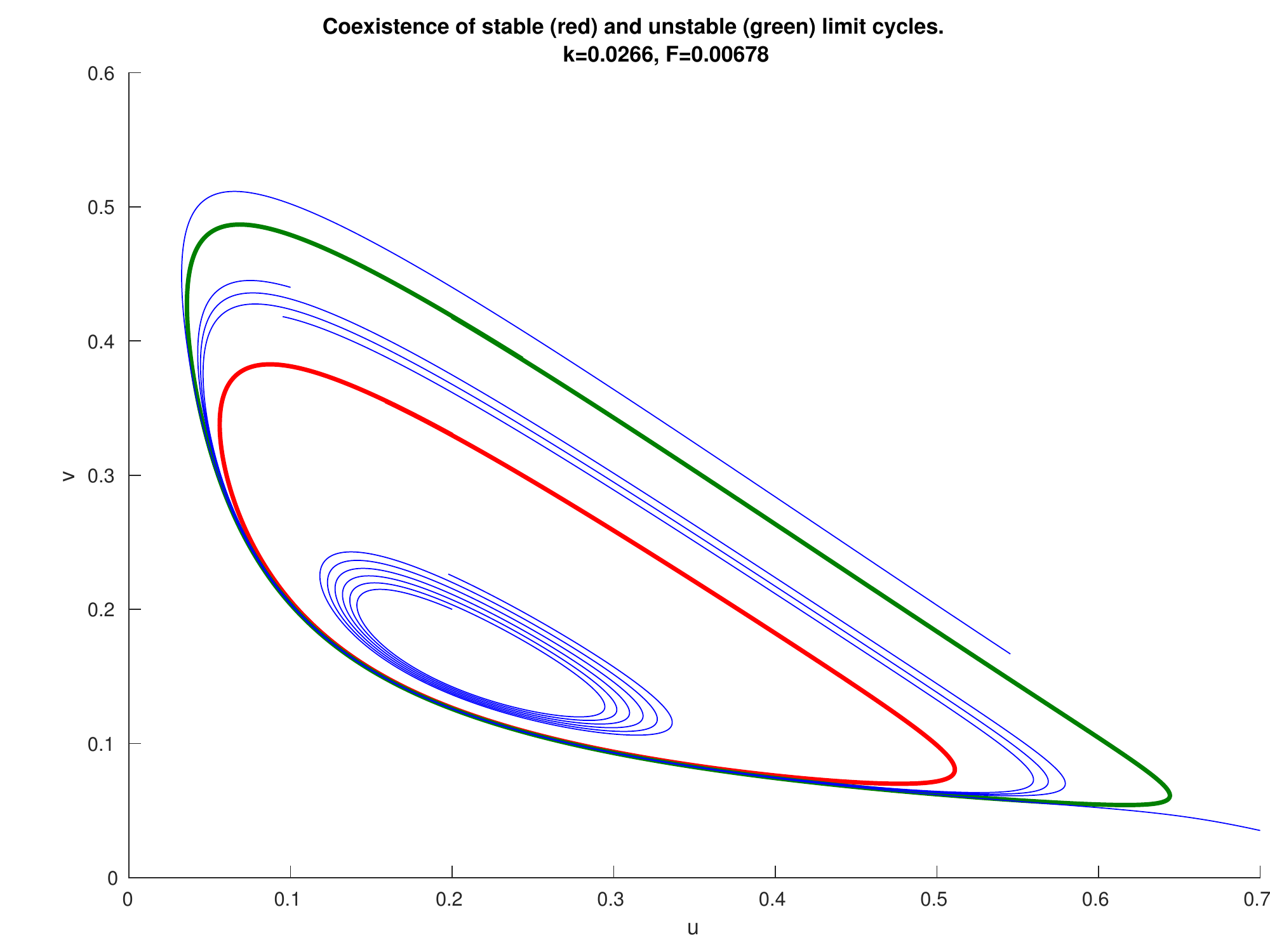}
\caption{Consequence of the existence of a point of bifurcation Bautin: A open region where coxisist two limit cycles. }
\end{figure}

\section{The global bifurcation map of the spatially homogeneous states}
\begin{figure}\centering
\includegraphics[width=\textwidth]{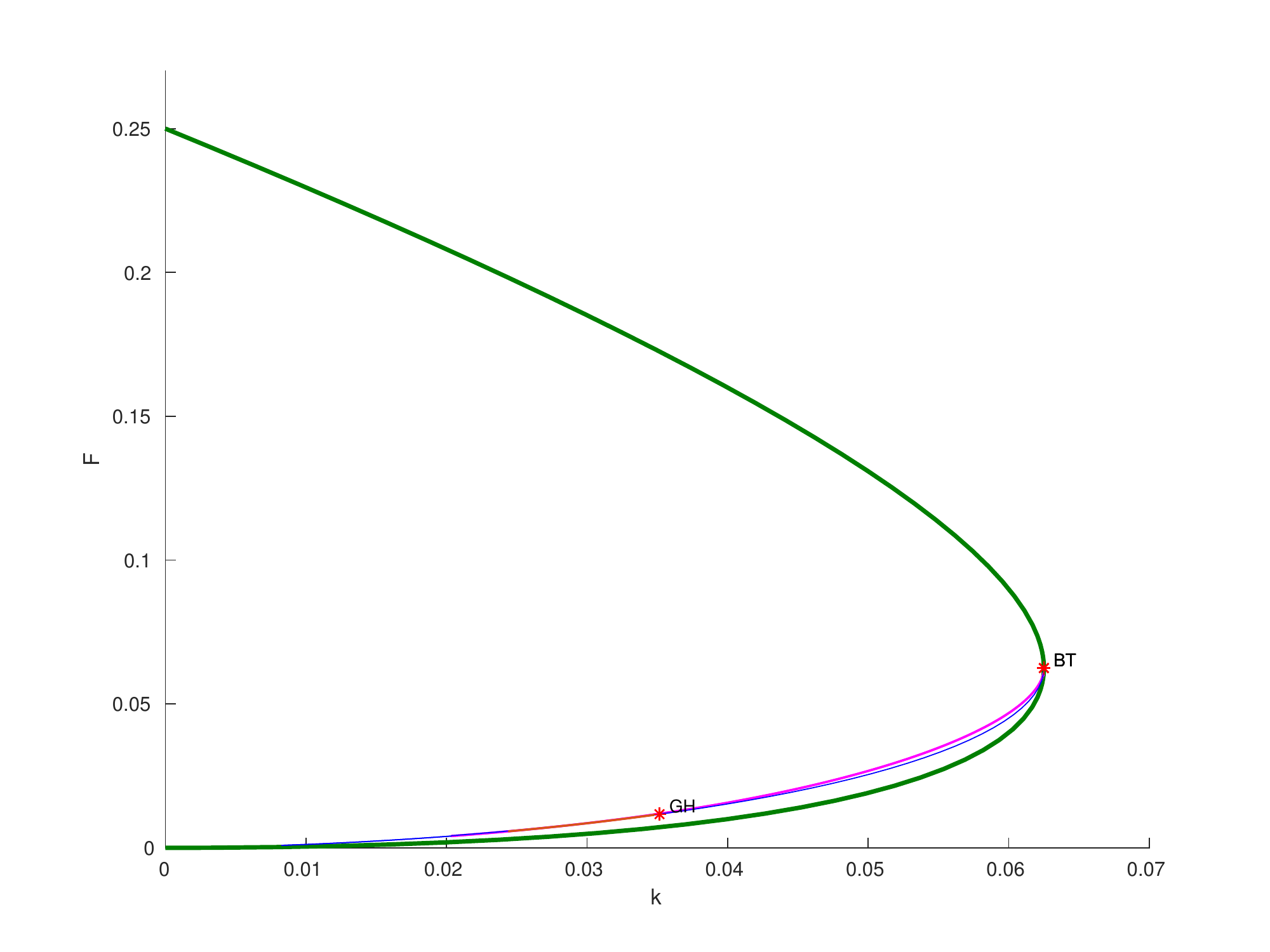}
\caption{The numerical continuation of the Hopf and homoclinic curves originating from the Takens-Bogdanov point BT \label{Numerical_continuation}}
\end{figure}

\begin{figure}\centering
\includegraphics[width=\textwidth]{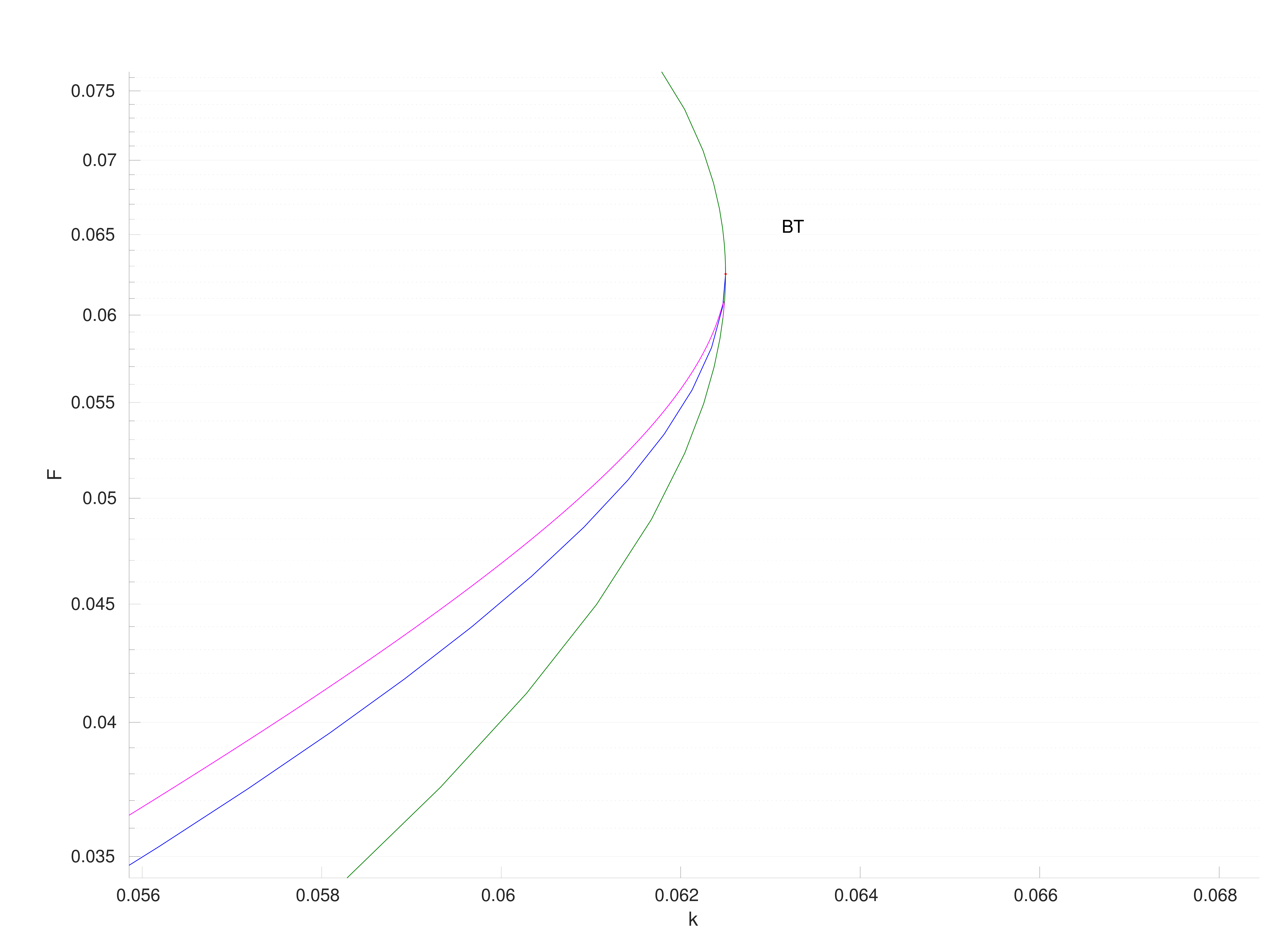}
\caption{Detail of figure \ref{Numerical_continuation} in a neighborhood of the BT point, showin the homoclinic, Hopf and saddle--nodes curves, from top to bottom (magenta, blue and green respectively) \label{global_map_zoom1} }
\end{figure}

According  to Theorem \ref{tTB}, there exist a local map of bifurcation curves in the parameter space, consisting of two components of saddle--nodes, and two curves of Hopf and a homoclinic points, both tangent to the saddle--node curve at BT. The numerical continuation of these curves, computed with MatCont,  are shown in Figure \ref{Numerical_continuation}. Detail of the relative position of the homoclinic, Hopf and saddle--node curve are shown in Figure~\ref{global_map_zoom1}.

Similarly, according to Theorem \ref{Bautin}, there exist  local curves of Hopf points and  a curve of limit point of cycles $T$, which is tangent to the Hopf curve. Actually, since we have shown that the second Lyapunov coefficient $\ell_2$ is positive at GH, see (\ref{L2}),  the local bifurcation map in a neigborhood of the point GH point is described in Figure~\ref{Bautin-bif}. Note that the numerical continuation of the curve $T$ is very near to the Hopf curve and can only be distinguished in detail as shown in Figure~\ref{global_map_zoom2}, right.

\begin{figure}\centering
\includegraphics[width=\textwidth]{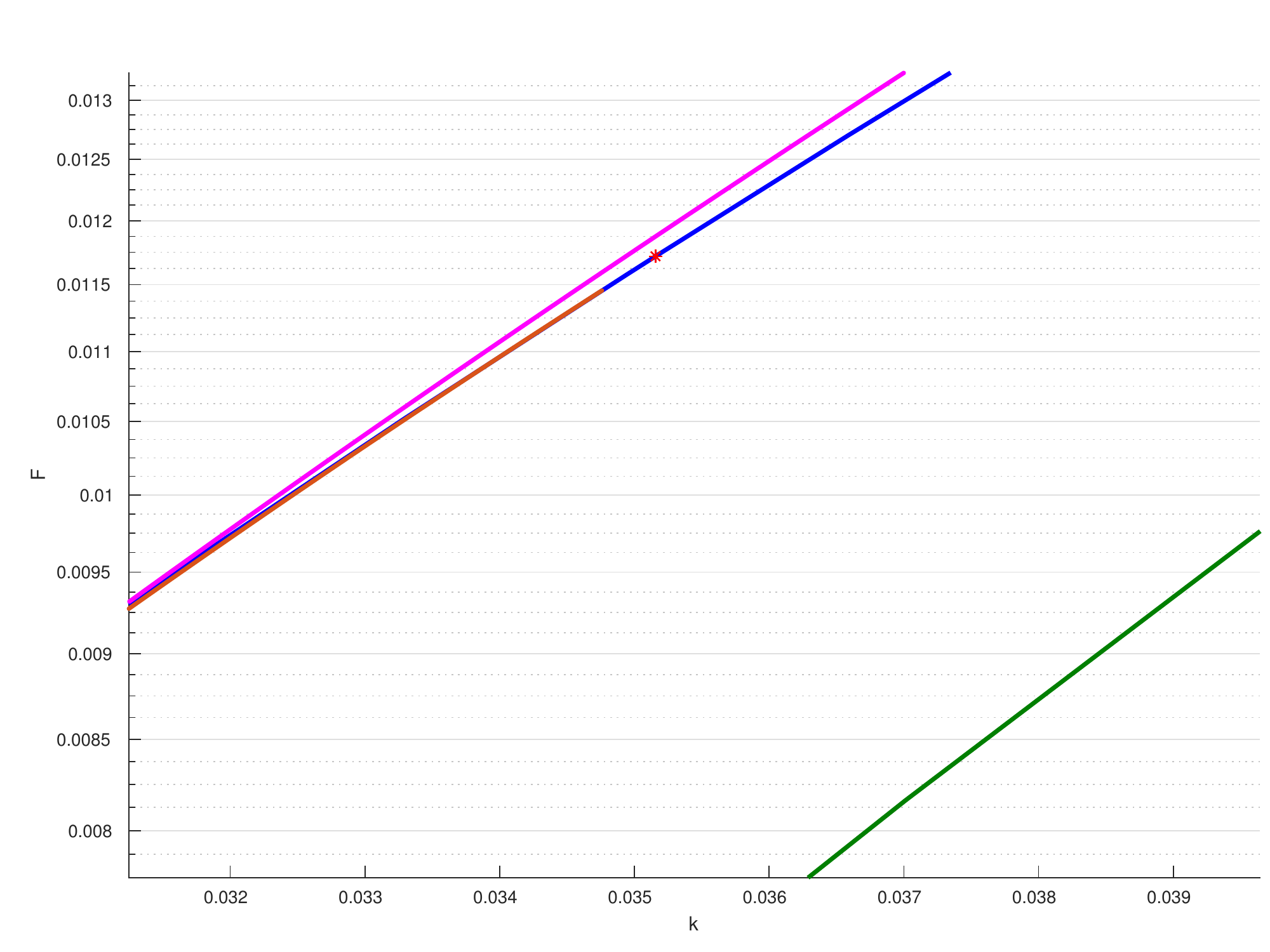}
\caption{Detail of figure \ref{Numerical_continuation} in a neighborhood of the Bautin point (red asterisc), the homoclinic and saddle--nodes curves are also shown with the same colors as In Figure \ref{global_map_zoom1}. The curve of limit point of cycles emerging from the Bautin point is shown as in organge color.  \label{intershomo-hopf} }
\end{figure}

\begin{figure}\centering
\includegraphics[width=0.48\textwidth, height=2.5in]{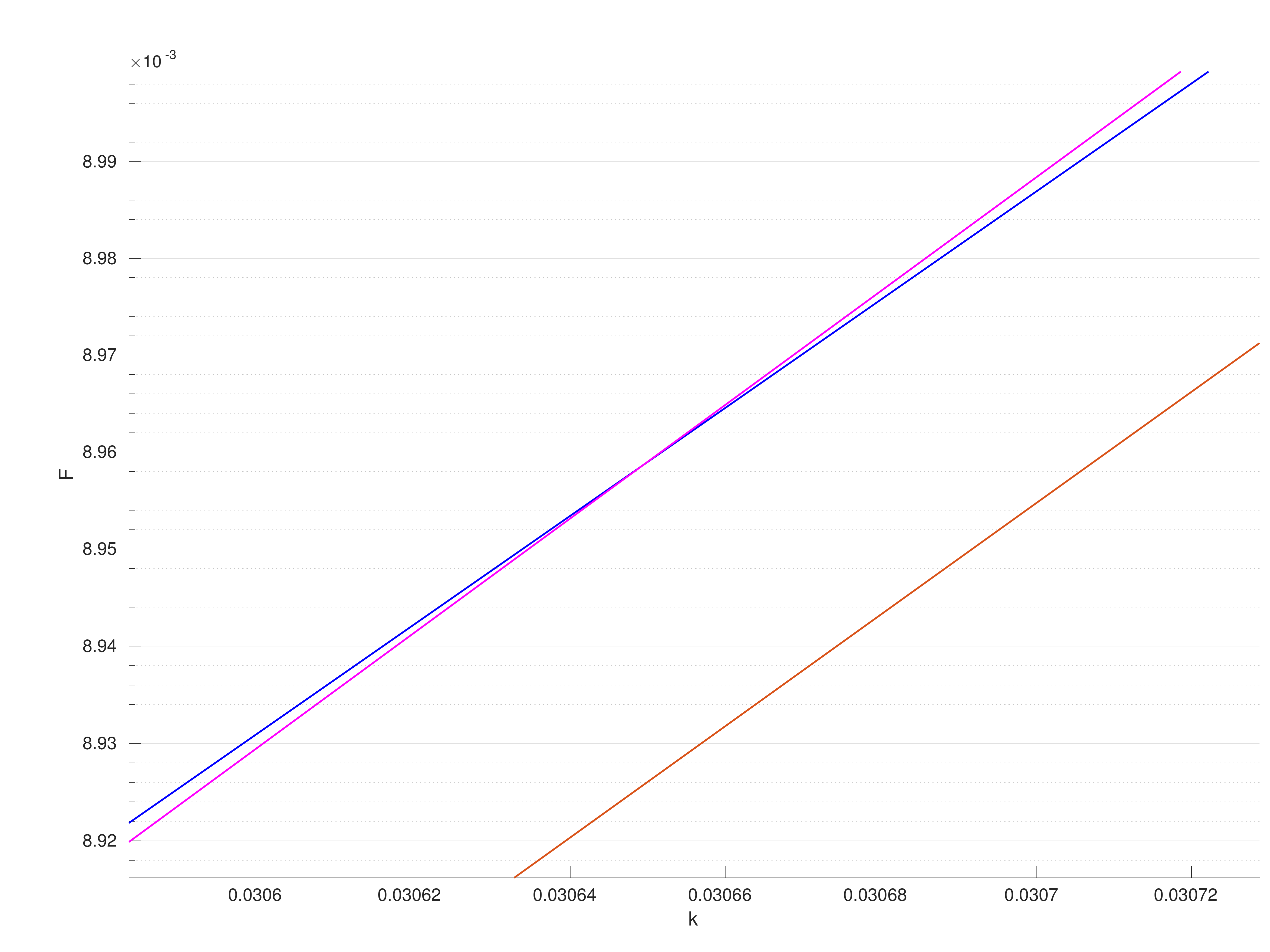}
\includegraphics[width=0.48\textwidth, height=2.5in]{zoom1}
\caption{The evolution of the homoclinic, Hopf and Bautin curve for different ranges of values of $k$ (note the horizontal scale).\label{global_map_zoom2} }
\end{figure}

\begin{figure}[t]\centering
\includegraphics[width=\textwidth]{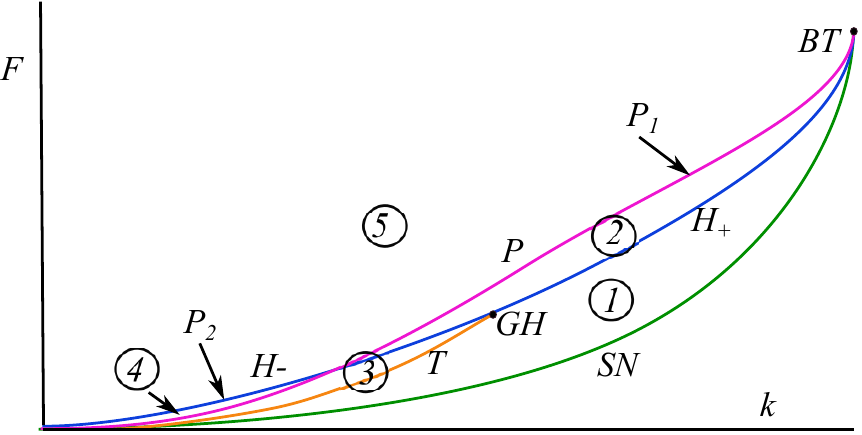}
\caption{Schematic diagram of the global bifurcation map.\label{Global-map-bif}}
\end{figure}

\begin{figure}[ht]
\begin{tabular}{|c|c|} \hline
\includegraphics[scale=.5]{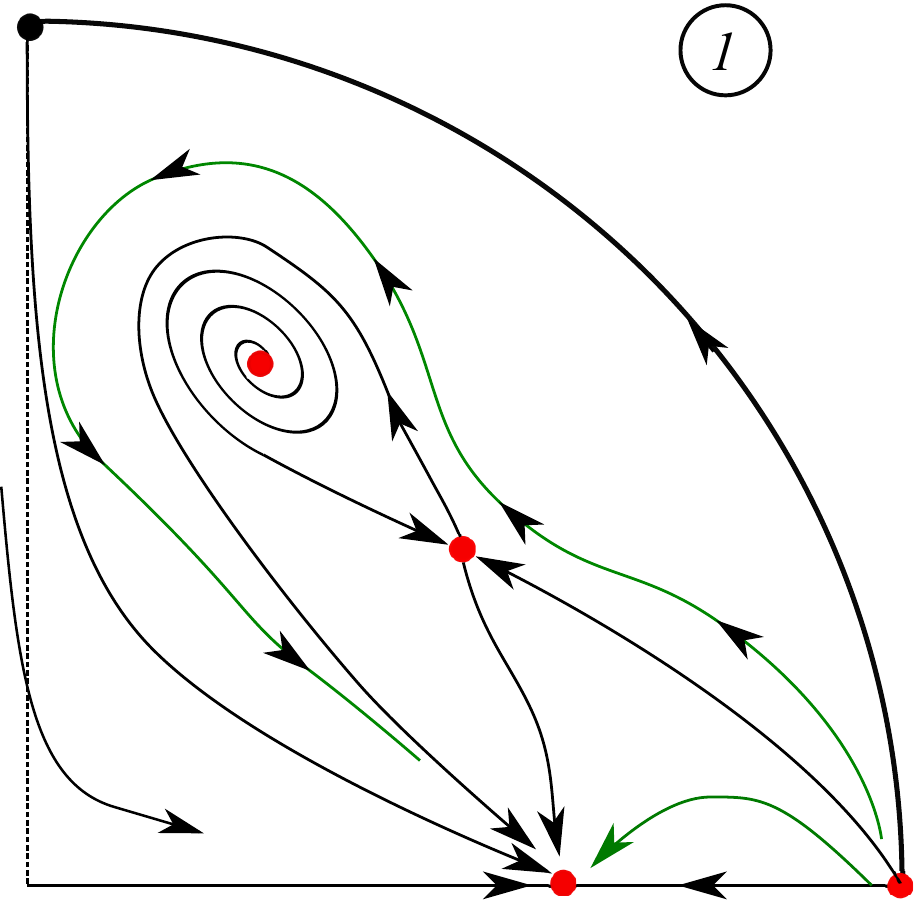}&
\includegraphics[scale=.5]{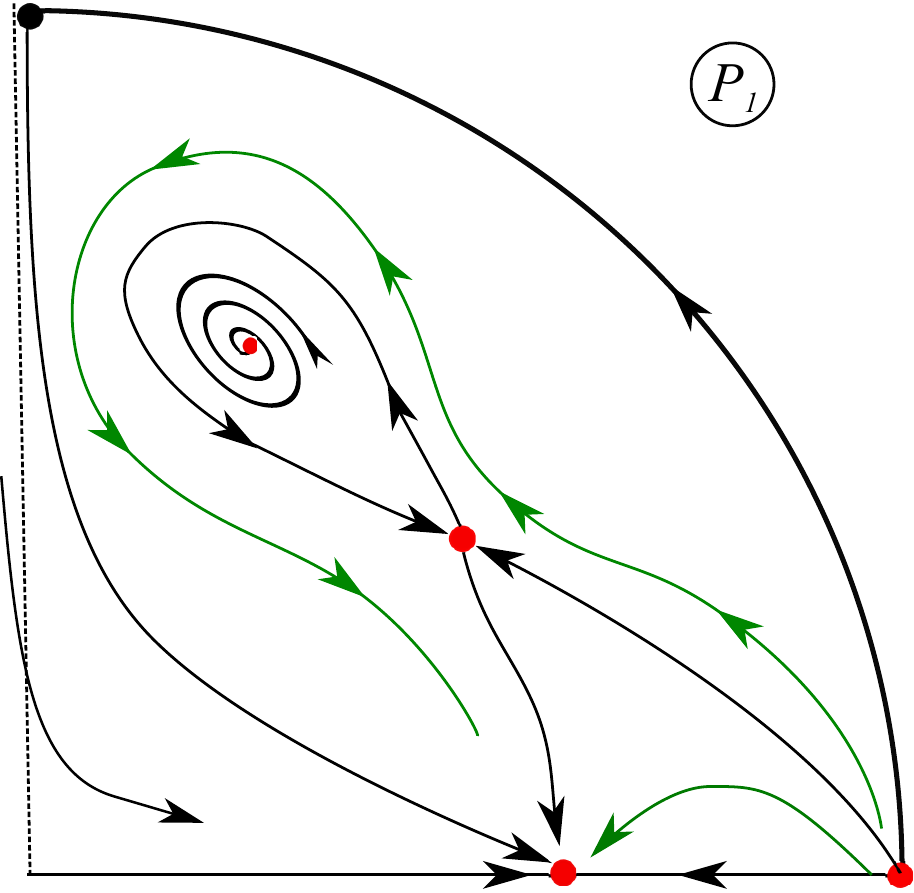}\\ \hline
\includegraphics[scale=.5]{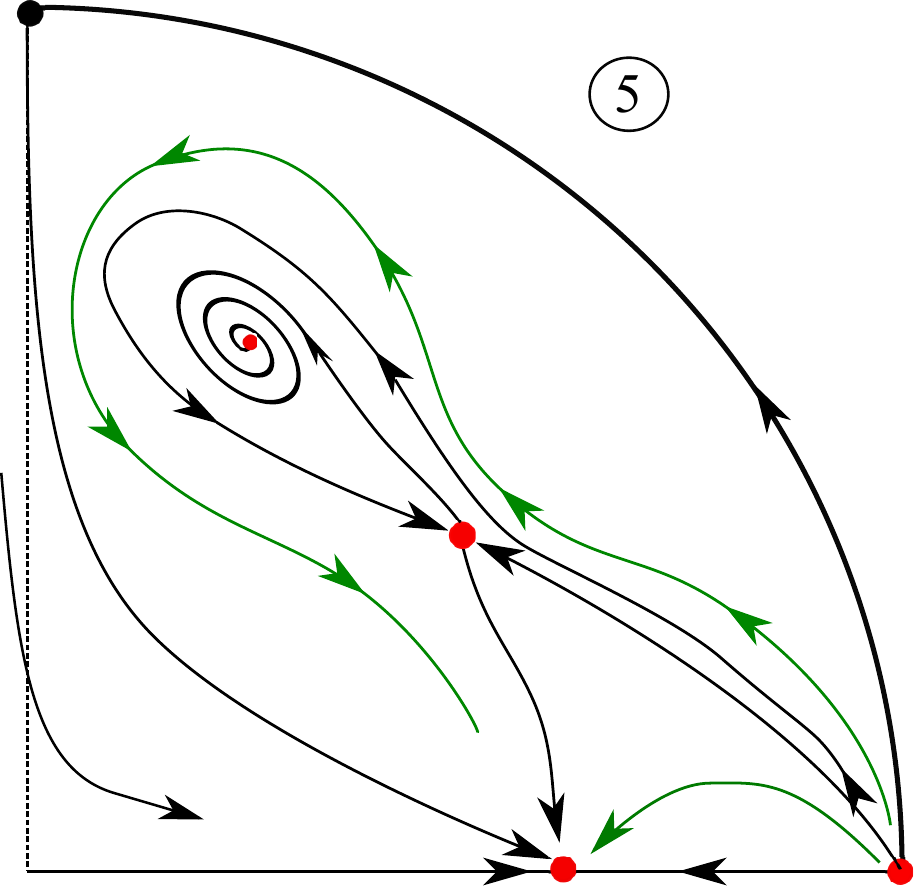}&
\includegraphics[scale=.5]{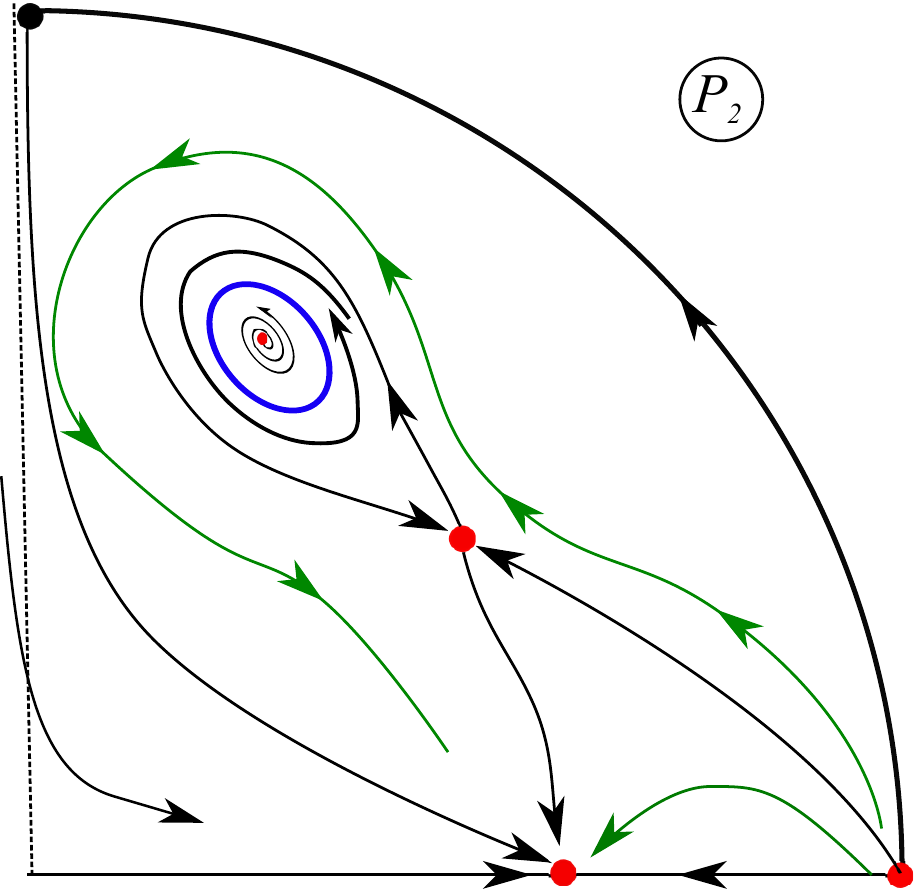}\\ \hline
\includegraphics[scale=.5]{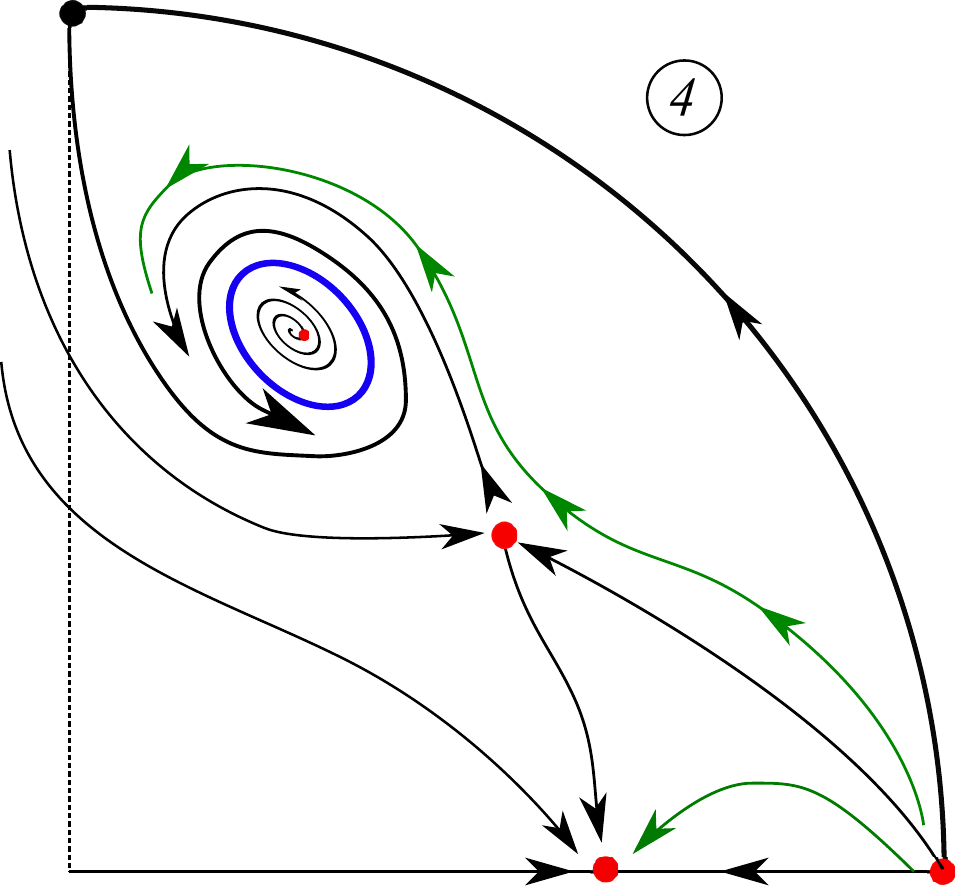}&
\includegraphics[scale=.5]{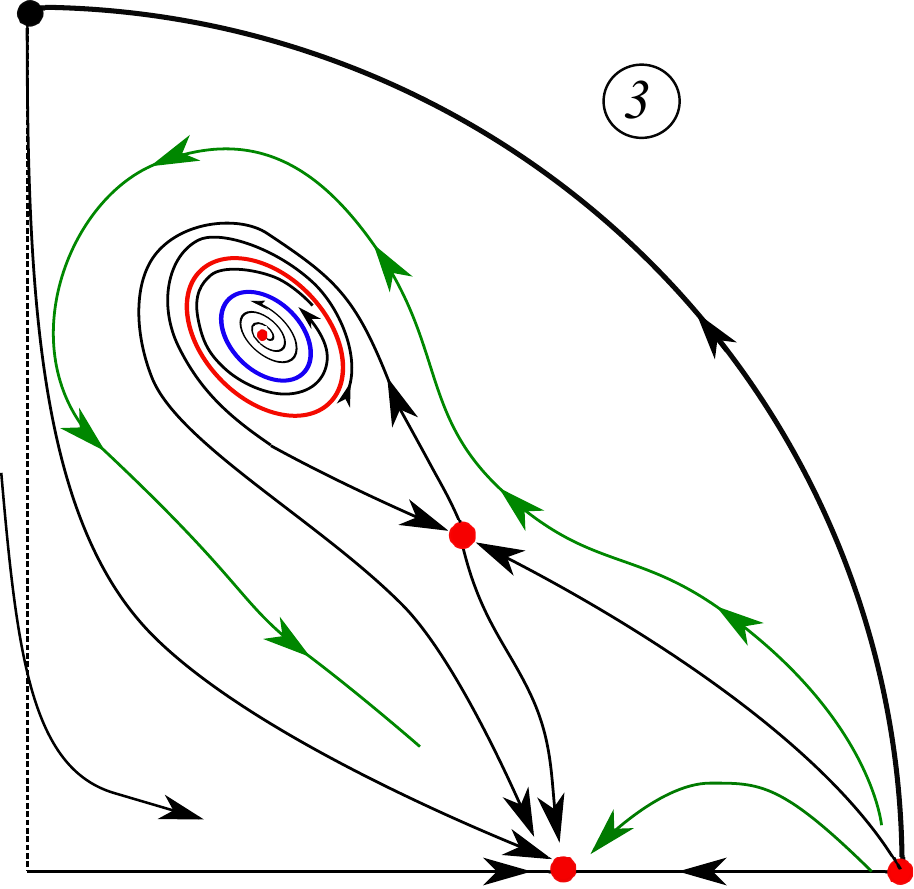}\\ \hline
\end{tabular}
\caption{Phase portrait of system (\ref{system}) for each or the regions numbered in Figure \ref{Global-map-bif}.\label{phase-portraits}}
\end{figure}

Figure~\ref{Global-map-bif} shows a schematic diagram of the bifurcation curves numerically continued according to Figures~\ref{Numerical_continuation}, \ref{global_map_zoom1} and \ref{global_map_zoom2}.  The saddle node curve (SN), the Hopf curves $H_{\pm}$, the homoclinic curve (P) and the Bautin curve (T) divide the interior of the saddle node curve into several regions, denoted 1 to 5 and especial points where a homoclinic bifurcation takes place, denoted by $P_1$ and $P_2$. The corresponding phase portraits are depicted in Figure~\ref{phase-portraits}.

The critical point $p_0$ lays along the $u$--axis and is an attractor for all values of the parameters. The critical point $p_{\pm}$ is a saddle and is marked as a red dot in the center of the phase plane. The critical point $p_{\mp}$ bifurcates changes stability as described in the global bifurcation diagram and generates families of limit cycles through Hopf bifurcations. The two limit cycles appear in the phase portrait numbered 3 in the Figure.

To complete the global phase portrait we complete the phase portrait by the Poincaré compactification (see details of this construction in \cite[p. 267]{Perko}). 
For each region in Figure~\ref{Global-map-bif}, Figure~\ref{phase-portraits} shows the phase portrait in the 2D Poincaré sphere where we add an invariant circle at infinity.
There are two  critical points at infinity marked as red and black dots which corresponds to asymptotic directions along $u=0$, $v=+\infty$ and $u=\infty$, $v=0$, respectively,   and a regular curve joining these points at infinity. The critical point corresponding to $u=0$, $v=+\infty$ is a degenerate saddle, for which there exists a one--dimensional unstable manifold which connects with  the trivial critical point $p_0$ for some value of the parameters and connects with the point $p_{\mp}$. Therefore the unstable manifold of the critical point at infinity $u=0$, $v=+\infty$ must connect with the saddle--point $p_{\pm}$ for some values of the parameters. 
We state this as a conjecture:
\bigskip

\noindent{\textbf{Conjecture.}} There exists a curve in the parameter space $k$-$F$ for which the unstable manifold of the critical point at infinity $u=0$, $v=+\infty$ connects with the stable manifold of the saddle point $p_{\pm}$.
\bigskip

\appendix
\section{Proof of Theorem \ref{tTB}. \label{TB}}
We start from the system (\ref{df}), $\dot x=f(x,\alpha)$ rewritten here for convenience
\begin{eqnarray*}
\dot x_1&=&-\left(\dfrac{(\alpha_2+\frac{1}{16})^2}{4(\alpha_1+\alpha_2+\frac{1}{8})^2}+\alpha_2+\frac{1}{16}\right)x_1-\dfrac{\alpha_2+\frac{1}{16}}{2(\alpha_2+\alpha_1+\frac{1}{8})}x_2\\
&-&\dfrac{\alpha_2+\frac{1}{16}}{\alpha_2+\alpha_1+\frac{1}{8}}x_1x_2-\dfrac{1}{2}{x_2}^2-x_1{x_2}^2, \\
\dot x_2&=&\dfrac{(\alpha_2+\frac{1}{16})^2}{4(\alpha_2+\alpha_1+\frac{1}{8})^2}x_1+\left(\dfrac{\alpha_2+\frac{1}{16}-2(\alpha_2+\alpha_1+\frac{1}{8})^2}{2(\alpha_2+\alpha_1+\frac{1}{8})}\right)x_2\\
&+&\dfrac{1}{2}x_2^2+\dfrac{\alpha_2+\frac{1}{16}}{\alpha_2+\alpha_1+\frac{1}{8}}x_1x_2+x_1x_2^2,
\end{eqnarray*}

Let $A_{0}=D_xf(0,0)$ and  $F(x)=Df(x,0)-A_{0} x$ its non--linear part,
$$A_0=\left(\begin{array}{cc}
-\frac{1}{8}&-\frac{1}{4}\\
\frac{1}{16}&\frac{1}{8}
\end{array}\right),\hspace{3mm}
F(x)=\left(\begin{array}{c}
-\frac{x_1 x_2}{2}-\frac{x_2^2}{2}-x_1x_2^2\\
\frac{x_1x_2}{2}+\frac{x_2^2}{2}+x_1 x_2^2
\end{array}\right).
$$

Observe that $\det(A_0)=0,\mbox{ tr}(A_{0})=0$ and  $A_{0}\neq 0$. This imply that there  exist linearly independent vectors $\{v_0,v_1\}$ such that $A_0 v_0=0$, $A_0v_1=v_0$, and a dual basis $\{w_0,w_1\}$ such that $A^{T}_{0}w_{1}=0$,   $A^{T}_{0}w_{0}=w_1$.  We can choose the bases such that
\begin{eqnarray}
\langle v_0,w_0\rangle=\langle v_1,w_1\rangle&=&1,\label{eq:v0w0}\\
\langle v_1,w_0\rangle=\langle v_0,w_1\rangle&=&0.\nonumber
\end{eqnarray}
where $\langle,\rangle$ denotes the standard inner product.
We choose
$$
\begin{array}{rclrrcl}
v_0 &=& (-2,1),  &v_1 &=& (0,8),\\
w_0 &=&(-1/2,0), & w_1&=&(1/16,1/8)
\end{array}
$$
Introduce coordinates $y_1$, $y_2$ relative to the basis $\{v_0,v_1\}$,
$$
x=y_1v_0+y_2v_1,
$$
which can be calculated as
\begin{eqnarray*}
y_1&=&\langle x,w_0\rangle,\\
y_2&=&\langle x,w_1\rangle.\\
\end{eqnarray*}
In this coordinates the system
$$\dot{ x}=A_{0}x+F(x),$$
transforms into
$$\dot{ y}=J y+G(y),$$
where
$$
J = B^{-1}A_{0}B = \left(\begin{array}{cc}  0 & 1 \\ 0 & 0 \end{array}\right)
$$
and  $B=(v_0|v_1)$. Explicitly,
$$\left(\begin{array}{c}
\dot y_1\\
\dot y_2
\end{array}\right)=\left(\begin{array}{cc}
0&1\\
0&0
\end{array}
\right)\left(\begin{array}{c}
y_1\\
y_2
\end{array}\right)+\left(\begin{array}{cc}
-\frac{y_1^2}{4}+16 y_2^2-y_1^3-16 y_2 y_1^2-64 y_2^2 y_1\\
-\frac{y_1^2}{32}+2 y_2^2-\frac{y_1^3}{8}-2 y_2 y_1^2-8 y_2^2 y_1
\end{array}
\right)
$$
where, according to \cite{Kuznetsov}, the second order terms $f(y)$ were obtained by the expression
$$
G(y)=\left(\begin{array}{c}
\langle F(y_1v_0+y_2v_1),w_0\rangle \\
\langle F(y_1v_0+y_2v_1),w_1\rangle
\end{array}\right)$$

The coordinates $(y_1,y_2)$ were obtained for the vector of parameters $\alpha=(0,0)$,  we now use these coordinates for the full parameterized vector field $\dot x=f(x,\alpha)$ for small $\alpha$, in the form
\begin{equation}\label{eq:fbeta}
\left(\begin{array}{c}
\dot y_1\\
\dot y_2
\end{array}\right)=\left(\begin{array}{c}
\langle f(y_1v_0+y_2v_1,\alpha),w_0\rangle \\
\langle f(y_1v_0+y_2v_1,\alpha),w_1\rangle
\end{array}
\right)\end{equation}
Using (\ref{df}) in the above expressions we obtain
\begin{eqnarray*}
\langle f(y_1v_0+y_2v_1,\alpha),w_0\rangle&=&-\frac{(16 \alpha_2+1)\left(4 \alpha_2^2+8 \alpha_2 \alpha_1+\alpha_2+4 \alpha_1^2\right)}{(8 \alpha_2+8 \alpha_1+1)^2} y_1\\
&+&\frac{(16 \alpha_2+1)}{8 \alpha_2+8 \alpha_1+1}y_2-\frac{(24 \alpha_2-8 \alpha_1+1)}{4 (8 \alpha_2+8 \alpha_1+1)}y_1^2\\
&+&\frac{32 (\alpha_1-\alpha_2)}{8 \alpha_2+8 \alpha_2+1}y_1 y_2 +16 y_2^2\\
&-&16 y_2 y_1^2-64 y_2^2 y_1
\end{eqnarray*}
and
\begin{eqnarray*}
\langle f(y_1v_0+y_2v_1,\alpha),w_1\rangle&=&-\frac{(16 \alpha_1+1)\left(4 \alpha_2^2+8 \alpha_2 \alpha_1+4 \alpha_1^2+\alpha_1\right)}{8 (8 \alpha_2+8 \alpha_1+1)^2} y_1 \\
&-&\frac{2\left(4 \alpha_2^2+8 \alpha_2 \alpha_1+4 \alpha_1^2+\alpha_1\right)}{8 \alpha_2+8 \alpha_1+1} y_2 \\
&-&\frac{(24 \alpha_2-8 \alpha_1+1)}{32 (8 \alpha_2+8 \alpha_1+1)}y_1^2+\frac{4(\alpha_1-\alpha_2)}{8 \alpha_2+8 \alpha_1+1} y_1 y_2 +2 y_2^2\\
&-&\frac{y_1^3}{8}-2 y_2 y_1^2-8 y_2^2 y_1.
\end{eqnarray*}

Let $a_{i,j}$ denote the coefficient of  $y^{i}_{1}y^{j}_{2}$ in $\langle f(y_1v_0+y_2v_1,\alpha),w_0\rangle$ and $b_{i,j}$ the coefficient of  $y^{i}_{1}y^{j}_{2}$ in  $\langle f(y_1v_0+y_2v_1,\alpha),w_1\rangle$. In order to verify the hypotheses of the Takens--Bogdanov theorem we need the following coefficients that can be obtained  as,
\begin{eqnarray*}
a_{20}(\alpha)&=&\dfrac{\partial^2}{\partial y^2_1}\left.\langle f(y_1v_0+y_2v_1,\alpha),w_0\rangle\right|_{y=0}=-2\frac{(24 \alpha_2-8 \alpha_1+1)}{4 (8 \alpha_2+8 \alpha_1+1)},\\
b_{20}(\alpha)&=&\dfrac{\partial^2}{\partial y^2_1}\left.\langle f(y_1v_0+y_2v_1,\alpha),w_1\rangle\right|_{y=0}=2\frac{(24 \alpha_2-8 \alpha_1+1)}{32 (8 \alpha_2+8 \alpha_1+1)},\\
b_{11}(\alpha)&=&\dfrac{\partial^2}{\partial y_1\partial y_2}\left.\langle f(y_1v_0+y_2v_1,\alpha),w_1\rangle\right|_{y=0}=\frac{4(\alpha_1-\alpha_2)}{8 \alpha_2+8 \alpha_2+1}.\\
\end{eqnarray*}
In particular for $\alpha=(0,0)$,
\begin{eqnarray*}
a_{20}(0)&=&-\frac{1}{2},\\
b_{20}(0)&=&\frac{1}{16},\\
b_{11}(0)&=&0.
\end{eqnarray*}
So that
\begin{eqnarray*}
a_{20}(0)+b_{11}(0) &=&-1/2\neq 0,\\
b_{20}(0) &=&1/16\neq 0.
\end{eqnarray*}
In this way conditions (BT.1) and (BT.2) of \cite[p. 321]{Kuznetsov} are satisfied.

The last condition to be verified is the non--degeneracy of the map
$$ (x,\alpha)\rightarrow \left(f(x,\alpha),Tr\left(\dfrac{\partial f(x,\alpha)}{\partial x}\right),\det\left(\dfrac{\partial f(x,\alpha)}{\partial x}\right)\right)^T.
 $$
We verify this condition in the original coordinates $(u,v)$ and vector of parameters $(k,F)$, since the change of variables $(u,v,k,F)\mapsto (x,\alpha)$ is just a shift, thus the mapping is
$$ (u,v,k,F)\mapsto \left(\begin{array}{c}
F (1-u)-u v^2\\
u v^2-v (F+k)\\
-2 F-k+2 u v-v^2\\
F^2+F k-2 F u v+F v^2+k v^2
\end{array}\right)
$$
and the Jacobian matrix is
$$\left(
\begin{array}{cccc}
 -v^2-F & -2uv& 0 & 1-u \\
 v^2 & -F-k+2 u v & -v &-v \\
 2v & 2(u-v) & -1 & -2 \\
-2Fv& -2Fu+2(F+k)v &F+v^2 &2F+k+v(-2u+v)\\
\end{array}
\right).$$
Evaluating at  $(u,v)=(1/2,1/4)$ and $(k,F)=(1/16,1/16)$ gives
$$
\left|
\begin{array}{cccc}
 -\frac{1}{8} & -\frac{1}{4} & 0 & \frac{1}{2} \\
 \frac{1}{16} & \frac{1}{8} & -\frac{1}{4} & -\frac{1}{4} \\
 \frac{1}{2} & \frac{1}{2} & -1 & -2 \\
 -\frac{1}{32} & 0 & \frac{1}{8} & 0 \\
\end{array}
\right|= -\frac{1}{512}.
$$
\hfill$\Box$

\section{Proof of Theorem \ref{Bautin} \label{B}}
In this section we give the details of the proof of Theorem~\ref{Bautin}.
In order to use Chow's formula let us write system (\ref{system}) as  its linear part plus high order terms $\dot x=A(\alpha)x+F(x,\alpha)$ where
$$
A=\left(
\begin{array}{cc}
a(\alpha)&b(\alpha)\\
c(\alpha)&d(\alpha)
\end{array}
\right), \qquad F(x,\alpha)=\left(
\begin{array}{c}
f(x,\alpha)\\
g(x,\alpha)
\end{array}
\right).
$$

By hypothesis the point $p_{\mp}$ satisfies  $\mbox{tr}(A)=0$ and $\beta^2_0=\det(A)>0$.
Then the first coefficient of Lyapunov $\ell_{1}^{CLW}$ is
\begin{eqnarray*}
& &\frac{b}{16\beta^4_0}\left\{\beta^2_0\left[b(f_{xxx}+g_{xxy})+2d(f_{xxy}+g_{xyy})-c(f_{xyy}+g_{yyy})\right]\right.\\
&-&bd(f^2_{xx}+f_{xx}g_{xy}-f_{xy}g_{xx}-g_{xx}g_{yy}-2g^2_{xy})\\
&-&cd(g^2_{yy}+g_{yy}f_{xy}-g_{xy}f_{yy}-f_{yy}f_{xx}-2f^2_{xy})\\
&+&b^2(f_{xx}g_{xx}+g_{xx}g_{xy})-c^2(f_{yy}g_{yy}+f_{xx}f_{yy})\\
&-&\left.(\beta_0^2+3d^2)(f_{xx}f_{xy}-g_{xy}g_{yy})\right\}
\end{eqnarray*}
evaluated at $p_{\mp}$. Then we perform the change of parameters
$$ \nu=F +\frac{\sqrt{k}}{2} (-1 + \sqrt{1 - 4 \sqrt{k}} + 2 \sqrt{k}) $$
so that $\nu=0$ corresponds to the Hopf curve. Then using the above formula we obtain the first Lyapunov coefficient along the Hopf curve,  up to a non--zero factor,
\begin{eqnarray*}
\ell_{1}^{CLW}&=&
8 \left(94-17 \sqrt{1-4
   \sqrt{k}}\right)
   k^{7/2}+\left(3816-1596
   \sqrt{1-4 \sqrt{k}}\right)
   k^{5/2}\\
   &+&\left(1416-924
   \sqrt{1-4 \sqrt{k}}\right)
   k^{3/2}+\left(4 \sqrt{1-4
   \sqrt{k}}-66\right) k^4\\
   &+&252
   \left(3 \sqrt{1-4
   \sqrt{k}}-10\right) k^3+110
   \left(15 \sqrt{1-4
   \sqrt{k}}-28\right)
   k^2\\&+&\left(286 \sqrt{1-4
   \sqrt{k}}-372\right)
   k+\left(52-46 \sqrt{1-4
   \sqrt{k}}\right) \sqrt{k}+3
   \left(\sqrt{1-4
   \sqrt{k}}-1\right),
\end{eqnarray*}
The zeros of this equation are given by the solution of the system of polynomial equations
\begin{eqnarray*}
Q_1&=&x^8 (4 y-66)+8 x^7 (94-17
   y)+252 x^6 (3 y-10)+x^5
   (3816-1596 y)\\
   &+&110 x^4 (15
   y-28)+x^3 (1416-924 y)+x^2
   (286 y-372)+x (52-46 y)\\
   &+&3(y-1)=0,\\
Q_2&=&-4 x-y^2+1=0,
\end{eqnarray*}
which is equivalent to change of variables
\begin{eqnarray*}
y&=&\sqrt{1-4\sqrt{k}},\\
x&=&\sqrt{k}.
\end{eqnarray*}
Taking the resultant with respect to $x$ we obtain the necessary condition
$$
\mbox{Res}(Q_1,Q_2,x)=2 (y-1)^{16} (2 y-1)=0.
$$
The root $y=1/2$ yields the solution
\begin{equation}\label{GHvalues}
k=\frac{9}{256},\qquad F=\frac{3}{256}.
\end{equation}
Substituting this values we get
\begin{equation}\label{Bau}
p_{\mp}=\left(\frac{1}{4},\frac{3}{16}\right).
\end{equation}

 For the second part, we use parameters $(k,F)$ and retain symbolic expression for the eigenbasis and dual eigenbasis, and make the substitution of the particular values of $k,F$ after computing the necesary derivatives of $\ell_1^{Kuz}(k,F)$ to prove the invertibility of the change of parameters
$$
(k,F)\mapsto (\beta_1,\beta_2)= (\mu(k,F),\ell_1^{Kuz}(k,F)),
$$
where $\mu(k,F) $ is the real part of the complex eigenvalue that vanishes at $GH$.
We get $-73728 \sqrt{2}$ for the Jacobian determinant at $GH$.

Another condition that has to be verified is that the second Lyapunov coefficient be different from zero. In this case wet get,
\begin{equation}\label{L2}
\ell_2^{Kuz}(GH)=10616832 \sqrt{2}.
\end{equation}
\hfill $\Box$

\medskip
\noindent{\bf \Large Acknowledgments}.
The second author would like to thank at Universidad Aut\'onoma Metropolitana-Iztapalapa and the Universidad Aut\'onoma de la Ciudad de M\'exico.

\end{document}